\documentclass[12pt]{article}
\usepackage{amsmath,amssymb,amsbsy,amsfonts,amsthm,latexsym,
                          amsopn,amstext,amsxtra,euscript,amscd}

\begin{document}

\newtheorem{theorem}{Theorem}
\newtheorem{lemma}[theorem]{Lemma}
\newtheorem{claim}[theorem]{Claim}
\newtheorem{cor}[theorem]{Corollary}
\newtheorem{prop}[theorem]{Proposition}
\newtheorem{definition}{Definition}
\newtheorem{question}[theorem]{Open Question}
\newtheorem{remark}{Remark}

\def\cP{\mathcal{P}}

\def\ep{\mathbf{e}_p}

\def\cT{\mathcal{T}}
\def\cU{\mathcal{U}}
\def\cV{\mathcal{V}}
\def\cA{\mathcal{A}}
\def\cX{\mathcal{X}}
\def\cY{\mathcal{Y}}
\def\cB{\mathcal{B}}
\def\cR{\mathcal{R}}
\def\cC{\mathcal{C}}
\def\cQ{\mathcal{Q}}
\def\cI{\mathcal{I}}

\title{Waring problem with the Ramanujan $\tau$-function}

\author{
{M.~Z.~Garaev, V. C. Garcia and S. V. Konyagin}}

\date{}%\today

\pagenumbering{arabic}

\maketitle

\begin{abstract}
Let $\tau(n)$ be the Ramanujan $\tau$-function. We prove that for
any integer $N$ the diophantine equation
$$
\sum_{i=1}^{74000}\tau(n_i)=N
$$
has a solution in positive integers $n_1, n_2,\ldots, n_{74000}$
satisfying the condition
$$
\max_{1\le i\le 74000}n_i\ll |N|^{2/11}+1.
$$
We also consider similar questions in the residue ring modulo a
large prime $p.$
\end{abstract}


\section{Introduction}

The Ramanujan  function $\tau(n)$ is defined by the expansion
$$
X\prod_{n=1}^{\infty}(1-X^n)^{24}=\sum_{n=1}^{\infty}\tau(n)X^n.
$$
It possesses many remarkable properties of arithmetical nature. It
is known that:
\begin{itemize}
\item{} $\tau(n)$ is an integer-valued multiplicative function, that
is, $\tau(nm)=\tau(n)\tau(m)$ if $\gcd(n,m)=1$;
\item{} For any integer $\alpha\ge 0$ and prime $q,$
$\tau(q^{\alpha+2})=\tau(q^{\alpha+1})\tau(q)-q^{11}\tau(q^\alpha).$
In particular, $\tau(q^2)=\tau^2(q)-q^{11}.$
\item{} $\tau(n)\equiv \sum\limits_{d|n}d^{11} \pmod {691}.$
\item{} $\tau(n)\equiv \sum\limits_{d|n}d^{11} \pmod {2^8}$ for $n$
odd.
\item{} $|\tau(q)|\le 2q^{11/2}$ for any prime $q$ and
$|\tau(n)|\ll n^{11/2+\varepsilon}$ for any positive integer $n$ and
any $\varepsilon>0$ (here and throughout the paper the constants
implicit in the Vinogradov symbols ``$\ll$" and ``$\gg$" may depend
only on $\varepsilon$); this has been proved by Deligne.
\end{itemize}
For these and other properties of $\tau(n),$ see for
example,~\cite{Del},~\cite{Iw},~\cite{Kobl},~\cite{Ser}. In
particular, one can derive from~\cite{Ser} that any residue class
modulo a prime number $p$ can be expressed as $\tau(n)\pmod p$ for
some positive integer $n$. Based on the deep sum-product estimate of
Bourgain, Katz and Tao~\cite{BKT}, and Vinogradov's double
exponential sum estimate, Shparlinski~\cite{Sh} established that the
values $\tau(n), n\le p^4,$ form a finite additive basis modulo $p$,
i.e., there exists an absolute integer constant $s\ge 1$ such that
any residue class modulo $p$ is representable in the form
$$
\tau(n_1)+\ldots+\tau(n_s)\pmod p
$$
with some positive integers $n_1,\ldots,n_s\le p^{4}.$ In the
present paper, introducing a new approach, we prove a result which
in a particular case reduces the number $4$ on the exponent of $p$
to the best possible $2/11.$

\begin{theorem}
\label{thm:RamWar} The set of values of $\tau(n)$ forms a finite
additive basis for the set of integers. Moreover, for any integer
$N$ the diophantine equation
$$
\sum_{i=1}^{74000}\tau(n_i)=N
$$
has a solution in positive integers $n_1, n_2,\ldots, n_{74000}$
satisfying the condition
$$
\max_{1\le i\le 74000}n_i\ll |N|^{2/11}+1.
$$
\end{theorem}

We remark that the quantity $74000$ comes from $36\times2050+200.$
Here the number $2050$ (as well as the number 200) can be reduced if
one uses the recent developments on the Waring-Goldbach problem,
see~\cite{ArCh, KT} and therein references. This however still gives
a big number of summands. In this connection we prove the following
results, where we assume that $p$ is a large prime number.

\begin{theorem}
\label{thm:Ram32} For any integer $\lambda$ the congruence
$$
\sum_{i=1}^{16}\tau(n_i)-\sum_{i=1}^{16}\tau(m_i)\equiv\lambda\pmod
p
$$
holds for some positive integers $n_1, m_1,\ldots, n_{16}, m_{16}$
with
$$
\max_{1\le i\le 16}\{n_i, m_i\}\ll p^2\log^4p
$$
and $\gcd(n_im_i, 23!)=1$.
\end{theorem}

Using Theorem~\ref{thm:Ram32}, we show that any residue class modulo
$p$ is representable in the form
$$
\sum_{i=1}^{96}\tau(n_i)\pmod p
$$
with $\max\limits_{1\le i\le 96}n_i\ll p^2\log^4p.$ In particular,
for some positive constant $C$ the set
$$
\{\tau(n)\pmod p: n\le Cp^2\log^4p\}
$$
forms a finite additive basis for the residue ring $\mathbb{Z}_p$ of
order at most $96,$ see~\cite{Kar} for the definition.

\begin{theorem}
\label{thm:Ram16} For any integer $\lambda$ and any $\varepsilon>0,$
the congruence
$$
\sum_{i=1}^{16}\tau(n_i)\equiv\lambda\pmod p
$$
is solvable in positive integers $n_1, \ldots, n_{16}$ with
$$
\max_{1\le i\le 16}n_i\ll p^{3+\varepsilon}.
$$
\end{theorem}
In particular, for any sufficiently large prime $p$ the set
$$
\{\tau(n)\pmod p: n\le p^{3+\varepsilon}\}
$$ is a basis of $\mathbb{Z}_p$ of order at most $16.$

\bigskip

Throughout the text the letters $q, q_1, q_2,\ldots,$ are used to
denote prime numbers. For a given set $\cA,$ $|\cA|$ denotes its
cardinality.

{\bf Acknowledgement}. The third author was supported by the INTAS
grant 03-51-5070.

\section{Preliminary statements}

First, we recall the following consequence of the classical result
of Hua Loo-Keng~\cite{Hua}.

\begin{lemma}
\label{lem:hua} Let $s_0$ be a fixed integer $\ge 2049$ and let $J$
denote the number of solutions of the Waring-Goldbach equation
$$
\sum_{i=1}^{s_0}q_i^{11}=N
$$
in primes $q_1, \ldots, q_{s_0}$ with $q_i>23$ for all $1\le i\le
s_0.$ There exist positive constants $c_1=c_1(s_0)$ and
$c_2=c_2(s_0)$ such that for all sufficiently large integer $N,
N\equiv s_0\pmod 2,$ the following bounds hold:
$$
c_1\frac{N^{s_0/11-1}}{(\log N)^{s_0}}\le J\le
c_2\frac{N^{s_0/11-1}}{(\log N)^{s_0}}.
$$
\end{lemma}

\bigskip

Next, we require the following result of Glibichuk~\cite{Gl}.
\begin{lemma}
\label{lem:Glib} If $\cX, \cY\in\mathbb{Z}_p$ with $|\cX||\cY|>2p,$
then
$$
\left\{\sum_{i=1}^{8}x_iy_i: \quad x_i\in \cX, \,
y_i\in\cY\right\}=\mathbb{Z}_p.
$$
\end{lemma}

We will also need some special computational results concerning the
values of $\tau(n).$ Due to the aforementioned result of Deligne
$$
|\tau(n)|\ll n^{11/2+\varepsilon},
$$
the following $N^6$ numbers
$$
\sum_{i=1}^6\tau(a_i), \quad 1\le a_1,\ldots,a_6\le N
$$
are all of the size $O(N^{11/2+\varepsilon}).$ Thus, in average,
each number is representable by  the sum of six values of $\tau(n)$
many times. It is natural to expect that zero can also be expressed
as a sum of six values of $\tau(n)$. With this in mind we search six
positive integers $a_1, \ldots, a_6$ satisfying
$$
\sum_{i=1}^6\tau(a_i)=0.
$$
There are many formulas which connect $\tau(n)$ with the function
$$
\sigma_s(n)=\sum_{d|n}d^s.
$$
It is known, for example, that
$$
\tau(n)=\frac{65}{756}\sigma_{11}(n)+\frac{691}{756}\sigma_5(n)-
\frac{691}{3}\sum_{k=1}^{n-1}\sigma_5(k)\sigma_5(n-k).
$$
Another formula  looks like
$$
\tau(n)=n^4\sigma_0(n)-24\sum_{k=1}^{n-1}(35k^4-52k^3n+18k^2n^2)\sigma_0(k)\sigma_0(n-k),
$$
see~\cite{Nie}. The formulas of the above type are useful for
numerical computations of $\tau(n).$ In particular, one can extract
that
$$
\tau(12)=-370944, \,\tau(27)=-73279080,\, \tau(55)=2582175960,
$$
$$
\tau(69)=4698104544,\, \tau(90)=13173496560,\,
\tau(105)=-20380127040.
$$
Now we have
\begin{equation}
\label{eqn:6tau}
\tau(12)+\tau(27)+\tau(55)+\tau(69)+\tau(90)+\tau(105)=0.
\end{equation}
Therefore, assuming that  Theorem~\ref{thm:Ram32} is proved, we
multiply the set equality
$$
\Bigl\{\sum_{i=1}^{16}\tau(n_i)-\sum_{i=1}^{16}\tau(m_i)\pmod p:
m_i,n_i\le Cp^2\log^4p, (m_in_i, 23!)=1\Bigr\} =\mathbb{Z}_p
$$
by
$$
\tau(12)=-\tau(27)-\tau(55)-\tau(69)-\tau(90)-\tau(105).
$$
Then using the multiplicative property of $\tau(n),$ we conclude
that the set
$$
\{\tau(n)\pmod p: n\le Cp^2\log^4p\}
$$
forms an additive basis of $\mathbb{Z}_p$ of order at most $96.$

It is also useful to note that
$$
\tau(6)=-6048, \,\tau(14)=401856,\, \tau(29)=128406630,
\,\tau(41)=308120442,
$$
$$
\tau(42)=101267712, \, \tau(44)=-786948864,\, \tau(48)=248758272.
$$
Thus, we have a representation of zero by the sum of seven values of
$\tau(n)$:
\begin{equation}
\label{eqn:7tau}
\tau(6)+\tau(14)+\tau(29)+\tau(41)+\tau(42)+\tau(44)+\tau(48)=0.
\end{equation}

\section{Proof of Theorem~\ref{thm:RamWar}}

Let $M$ be a large even integer parameter. Define the set
$$
\cQ=\{q:\,23<q\leq M^{1/11}\}.
$$
We call a subset $\cQ'\subset\cQ$  {\it admissible} if the equation
$$
\sum_{i=1}^6 \tau(q'_i)=\sum_{i=7}^{12}\tau(q'_i)
$$
has no solutions in  $q'_1,\dots,q'_{12}\in \cQ'$ satisfying
$$
q'_1<\ldots<q'_6;\quad q'_7<\ldots<q'_{12};\quad (q'_1,\ldots,
q'_6)\not=(q'_7,\ldots,q'_{12}).
$$
Using different properties of $\tau(n)$ it is easy to check that
there are admissible subsets with $12$ elements. This can be derived
if one combines the above mentioned congruences
$$
\tau(q)\equiv 1+q^{11}\pmod {691},\quad \tau(q)\equiv 1+q^{11}\pmod
{2^8},
$$
where $q$ is an odd prime, with chinese remainder theorem and prime
number theorem for arithmetical progressions to show that for any
$j, 1\le j\le 12,$ one can find a sufficiently large prime $\ell_j$
satisfying
$$
\tau(\ell_j)\equiv 2^j\pmod{8\times 691}.
$$

Next, among all admissible subsets we take $\cQ'$ to be one with the
biggest cardinality. If there are several such subsets, we take
$\cQ'$ to be one of them.

In particular, $|\cQ'|\ge 12$ and all the sums of the type
$$
\sum_{i=1}^6 \tau(q'_i): \quad q'_1<\dots<q'_6, \quad
q'_1,\dots,q'_6\in \cQ',
$$
are distinct. By Deligne's estimate for $\tau(q)$ and pigeonhole
principle we have
$$|\cQ'|^6\ll
(M^{1/11})^{11/2},
$$
whence,
\begin{equation}
\label{eqn:setminus est} |\cQ'|\ll M^{1/11-1/132}.
\end{equation}
Given $q\in \cQ\setminus\cQ',$ consider the set $\cQ'\cup \{q\}$.
From the maximality of $|\cQ'|$ we know that there exist
$q'_1,\ldots, q'_{12}\in\cQ'\cup \{q\}$ such that
\begin{equation}
\label{eqn:6=6 with q}
\sum_{i=1}^{6}\tau(q'_i)=\sum_{i=7}^{12}\tau(q'_i)
\end{equation}
and
\begin{equation}
\label{eqn:6;6 with q} q'_1<\ldots<q'_6; \quad q'_7<\ldots<q'_{12};
\quad (q'_1,\ldots, q'_6)\not=(q'_7,\ldots,q'_{12}).
\end{equation}
From the definition of $\cQ'$ we derive that
$$
q\in\{q'_1,\ldots, q'_{12}\}.
$$
Besides, due to~\eqref{eqn:6;6 with q}, $q$ occurs in the sequence
$q'_1, \ldots, q'_{12}$ at most two times. If $q$ occurs there
twice, then it occurs exactly one time in the sequence
$$
q'_1, \ldots, q'_6
$$
and exactly one time in the sequence
$$
q'_7, \ldots, q'_{12}.
$$
Therefore, we can cancel  both sides of~\eqref{eqn:6=6 with q} by
$\tau(q)$ and enumerating the remaining numbers $q_i'$, we obtain
that there exist
\begin{equation}
\label{eqn: 10 in Q'} q'_1, \ldots, q'_{5}\in\cQ';\quad q'_7,\ldots,
q'_{11}\in\cQ'
\end{equation}
such that
$$
\sum_{i=1}^{5}\tau(q'_i)=\sum_{i=7}^{11}\tau(q'_i)
$$
and
\begin{equation}
\label{eqn:Q enumarate} q'_1<\ldots<q'_5; \quad q'_7<\ldots<q'_{11};
\quad (q'_1,\ldots, q'_5)\not=(q'_7,\ldots,q'_{11}).
\end{equation}
Since $|\cQ'|>10,$ there exists
\begin{equation}
\label{eqn:exist q' in Q'} q'\in \cQ'\setminus\{q_1',\ldots,q'_{5},
q'_7,\ldots, q'_{11}\}.
\end{equation}
Then we have
$$
\sum_{i=1}^{5}\tau(q_i')+\tau(q')=\sum_{i=7}^{12}\tau(q_i')+\tau(q').
$$
This, in view of~\eqref{eqn: 10 in Q'},~\eqref{eqn:Q enumarate}
and~\eqref{eqn:exist q' in Q'}, contradicts the definition of
$\cQ'$.

Therefore, there is only the possibility that $q$ occurs in the
sequence
$$
q'_1, \ldots, q'_{12}
$$
exactly one time. Thus, in view of~\eqref{eqn:6=6 with q}, there
exist $\tilde q_1, \ldots,\tilde q_{11}\in \cQ'$ such that
$$
\sum_{i=1}^6\tau(\tilde q_i)=\sum_{i=7}^{11}\tau(\tilde
q_i)+\tau(q).
$$
Hence, for any $q\in\cQ\setminus\cQ'$ we get
\begin{eqnarray*}
\sum_{i=1}^6 \tau(\tilde q_iq)-\sum_{i=7}^{11}\tau(\tilde q_iq)-\tau(q^2)=\\
\tau(q)\sum_{i=1}^6\tau(\tilde
q_i)-\tau(q)\sum_{i=7}^{11}\tau(\tilde q_i)
-\tau(q^2)=\\
\tau^2(q)-\tau(q^2)=q^{11}.
\end{eqnarray*}

Thus we have proved that for any element $q\in\cQ\setminus\cQ'$
there exist elements $\tilde q_1,\ldots,\tilde q_{11}\in \cQ'$ such
that
\begin{equation}
\label{eqn:q11=6-6} q^{11}=\sum_{i=1}^6 \tau(\tilde
q_iq)-\sum_{i=7}^{11}\tau(\tilde q_iq)-\tau(q^2).
\end{equation}

Our next aim is to prove the solubility of the Waring-Goldbach
equation
\begin{equation}
\label{eqn:WQ} \sum_{j=1}^{2050}q_j^{11}= M
\end{equation}
in primes $q_1,\ldots,q_{2050}\in\cQ\setminus\cQ'.$  First of all,
from~\eqref{eqn:WQ} we have $q_j\le M^{1/11}.$ Next, for the number
of solutions of~\eqref{eqn:WQ} with $q_j\in\cQ$ we have, according
to Lemma~\ref{lem:hua} with $s_0=2050,$ the lower bound
\begin{equation}
\label{eqn:Number Sol Q} \ge c_1\frac{M^{2050/11-1}}{(\log
M)^{2050}}.
\end{equation}
For the number of solutions of~\eqref{eqn:WQ} with at least one
$q_{j_0}\in\cQ'$ we have, according to Lemma~\ref{lem:hua} with
$s_0=2049$ and~\eqref{eqn:setminus est}, the upper bound
\begin{equation}
\label{eqn:Number Sol Q'} \le
2050c_2|\cQ'|\frac{M^{2049/11-1}}{(\log M)^{2049}}\ll
\frac{M^{2050/11-1-1/132}}{(\log M)^{2049}}.
\end{equation}
Thus,  $\eqref{eqn:Number Sol Q}>\eqref{eqn:Number Sol Q'}.$
Therefore, \eqref{eqn:WQ} is solvable in $q_j\in\cQ\setminus\cQ'.$
We fix one of such solutions $(q_1, \ldots, q_{2050})$. To each
$q_j, \, 1\le j\le 2050,$ we apply~\eqref{eqn:q11=6-6} with $q=q_j$
and then perform the summation over $1\le j\le 2050.$ Since $\tilde
q_iq\le M^{2/11},$ we obtain
$$
M=\sum_{i=1}^{6\times2050}\tau(n_i)-\sum_{i=1}^{6\times2050}\tau(m_i),
$$
where
$$
\max_{1\le i\le 6\times 2050}\{n_i, m_i\}\le M^{2/11}, \quad
\gcd(n_im_i,23!)=1.
$$
Our assumption on $M$ is that it is a large even integer. Clearly,
we can multiply the above equality by $-1$ to have the same type of
representation for $-M$ as well. Furthermore, we can remove the
parity condition on $M$ by extracting one element $\tau(1)=1$ or
$\tau(29)\equiv 0\pmod 2$. Thus, any integer $M$ with sufficiently
large $|M|$ can be represented in the form
\begin{equation}
\label{eqn:M=tau-tau}
M=\sum_{i=1}^{6\times2050}\tau(n_i)-\sum_{i=1}^{6\times2050+1}\tau(m_i),
\end{equation}
where
$$
\max_{1\le i\le 6\times 2050+1}\{n_i, m_i\}\le |M|^{2/11}+1, \quad
\gcd(n_im_i,23!)=1.
$$
Recall that (see~\eqref{eqn:6tau})
$$
-\tau(12)=370944=\tau(27)+\tau(55)+\tau(69)+\tau(90)+\tau(105).
$$
Therefore, multiplying~\eqref{eqn:M=tau-tau} by $-\tau(12)$ and
using the multiplicative property of $\tau(n),$ we obtain
\begin{equation}
\label{eqn:M=tau+tau}
370944M=\sum_{i=1}^{6\times6\times2050+1}\tau(n_i)
\end{equation}
with $\max\limits_{1\le i\le 6\times6\times2050+1}n_i\le
106|M|^{2/11}.$

Let us show that any integer $r, 0\le r< 370944,$ can be expressed
as a sum of say exactly $198$ numbers of the form $\tau(n), n\le
105$ (an extra effort would reduce $198$ to a much smaller constant,
which however do not essentially influence to our final result). To
this end, we recall that
$$
\tau(1)=1,\, \tau(2)=-24,\, \tau(3)=252,\, \tau(5)=4830,
\,\tau(8)=84480.
$$
Now if $0\le r<370944,$ then
$$
r=84480r_5+r'_4=\tau(8)r_5+r'_4
$$
for some integers $0\le r_5\le 4$ and $0\le r'_4< 84480.$ Next, any
such $r'_4$ can be expressed as
$$
r'_4=4830r_4+r'_3=\tau(5)r_4+r'_3,
$$
where $0\le r_4\le 17$ and $0\le r'_3<4830.$ Any such $r'_3$ can be
written in the form
$$
r'_3=252r_3-r'_2=\tau(3)r_3-r'_2
$$
where $0\le r_3\le 20$ and $0\le r'_2< 252.$ Any such $r'_2$ can be
written in the form
$$
r'_2=24r_2-r_1=-\tau(2)r_2-r_1,
$$
with $0\le r_2\le 11$ and $0\le r_1\le 23.$ Thus we have
$$
r=\tau(8)r_5+\tau(5)r_4+\tau(3)r_3+\tau(2)r_2+\tau(1)r_1.
$$
Therefore, to express such a given $r,$ at most
$$
r_5+\ldots+r_1\le 75
$$
number of summands of $\tau(n), n\le 10,$ are sufficient. On the
other hand, any integer greater than $29$ can be expressed in the
form $6x+7y$ with nonnegative integers $x,y.$ Thus, in order to have
a fixed number of summands for all $r$, we can use
constructions~\eqref{eqn:6tau} and~\eqref{eqn:7tau}. In particular,
any integer $r, 0\le r< 370944,$ can be expressed in the form
$$
\sum_{i=1}^{198}\tau(a_i)
$$
with positive integers $a_1,\ldots,a_{198}\le 105.$

Let now $N$ be an arbitrary integer with a sufficiently large modulo
$|N|.$ The above argument shows that for some positive integers
$a_1,\ldots,a_{198}\le 105$ we have
$$
N\equiv \sum_{i=1}^{198}\tau(a_i)\pmod{370944}.
$$
Therefore, using~\eqref{eqn:M=tau+tau}, we obtain that
$$
N=\sum_{i=1}^{198}\tau(a_i)+370944M=\sum_{i=1}^{6\times6\times2050+199}\tau(n_i)=
\sum_{i=1}^{73999}\tau(n_i),
$$
where
$$
\max\limits_{1\le i\le 73999}n_i\le 106|M|^{2/11}\le 15|N|^{2/11}.
$$

Thus, we have proved that there exists an absolute positive integer
constant $N_0$ such that for any integer $N$ with $|N|\ge N_0$ the
equation
$$
\sum_{i=1}^{73999}\tau(n_i)=N
$$
has a solution in positive integers $n_1,\ldots,n_{73999}\ll
|N|^{2/11}$.

Let now $N$ be an arbitrary integer. If $|N|>N_0$ then
$|N-\tau(1)|\ge N_0$ and therefore we can express the number
$N-\tau(1)$ as a sum of $73999$ values of $\tau(n)$ with $n\ll
|N|^{2/11}.$ Theorem~\ref{thm:RamWar} follows in this case.

If $|N|\le N_0$, then we take an integer constant $n_0$ such that
$|\tau(n_0)|>2N_0$. Then
$$
|N-\tau(n_0)|>N_0.
$$
Thus $N-\tau(n_0)$ can be expressed as a sum of $73999$ values of
$\tau(n), \, n\ll 1.$

Theorem~\ref{thm:RamWar} is proved.

\begin{remark}
\label{rem:1} One can easily see that the numbers $n_i$ constructed
in the proof satisfy the condition
$$\tau(n_i)\ll |N|.$$
\end{remark}

\begin{remark}
\label{rem:2} In our forthcoming paper we will prove that for any
integer $N$ with $|N|\ge 2$ the diophantine equation
$$
\sum_{i=1}^{148000}\tau(n_i)=N
$$
has a solution in positive integers $n_1, n_2,\ldots, n_{148000}$
satisfying the condition
$$
\max_{1\le i\le 148000}n_i\ll |N|^{2/11}e^{-c\log |N|/\log\log |N|},
$$
for some absolute constant $c>0.$ In view of Deligne's estimate
$\tau(n)\le n^{11/2}d(n)$, where $d(n)$ is the number of divisors of
$n$, this reflects the best possible bound for the size of the
variables $n_i$, apart from the value of the constant $c$.
\end{remark}

\section{Proof of Theorem~\ref{thm:Ram32}}

Let $C$ be a large positive constant. Consider the sets
$$
\cQ=\{q: 23< q\le Cp^{1/2}\log p\}
$$
and
$$
\cI=\{\tau(q)\pmod p: q\in \cQ\}.
$$
If $|\cI|>3\sqrt{p},$ then we can split $\cI$ into two subsets $\cX,
\cY$ with $|\cX||\cY|>2p.$ The result in this case follows from
Lemma~\ref{lem:Glib}.

Let now $|\cI|<3\sqrt{p}.$ Then
$$
\cQ=\bigcup_{i=1}^{|I|}\cA_i,
$$
where the sets $\cA_i$ are defined such that the condition $q',
q''\in \cA_i$ implies $\tau(q')\equiv\tau(q'')\pmod p.$

Clearly, for some $\cA'_i\subset\cA_i$ we have
$$
0\le |\cA_i|-|\cA'_i|\le 3, \quad |\cA'_i|\equiv 0\pmod 4.
$$
Denote
$$
\cQ'=\bigcup_{i=1}^{|I|}\cA'_i.
$$
By the prime number theorem we have $|\cQ|\ge Cp^{1/2}$ provided
that $p$ is large enough. Then
\begin{equation}
\label{eqn:est Q' and Q} |\cQ'|\ge |\cQ|-3|I|\ge |\cQ|-9p^{1/2}\ge
(C-9)p^{1/2}.
\end{equation}
Since the cardinality of each set $\cA'_i$ is even, we can produce
$|\cA'_i|/2$ pairs formed with different primes of the set $\cA'_i.$
Then there are totally
$$
\sum_{i=1}^{|I|}|\cA'_i|/2=|\cQ'|/2
$$
pairs $(q,q').$ We divide this set of pairs into two disjoint
subsets $J_1$ and $J_2$, so that $|J_1|=|J_2|=|\cQ'|/4.$ Now
consider the following two sets:
$$
\cX=\{\tau(qq')-\tau(q^2)\pmod p:\quad (q, q')\in J_1\}
$$
and
$$
\cY=\{\tau(qq')-\tau(q^2)\pmod p:\quad (q, q')\in J_2\}.
$$
Since
$$
\tau(qq')-\tau(q^2)=\tau(q)\tau(q')-\tau(q^2)\equiv
\tau^2(q)-\tau(q^2)\equiv q^{11}\pmod p,
$$
and $q^{11}\pmod p$ can take any value at most 11 times, we have
$$
|\cX|\ge |J_1|/11=|\cQ'|/44, \quad |\cY|\ge |J_2|/11\ge |\cQ'|/44.
$$
Therefore, if we choose $C=100$ say, then according
to~\eqref{eqn:est Q' and Q}, we obtain $|\cX||\cY|>2p.$ Applying
Lemma~\ref{lem:Glib}, we finish the proof of
Theorem~\ref{thm:Ram32}.

\section{Proof of Theorem~\ref{thm:Ram16}}

Consider the set of residue classes
$$
\cA'=\{\tau(q)\pmod p:\quad p/2<q\le p\}.
$$
From different properties of $\tau(n)$ it follows that $\cA'$
contains more than one element (apply, for example, the above
mentioned congruence modulo $691$ and the prime number theorem for
arithmetical progressions). For a given $a'\in\cA',$ let $I(a')$ be
the number of solutions of the congruence
$$
\tau(q)\equiv a'\pmod p, \quad p/2<q\le p.
$$
From prime number theorem,
$$
\sum_{a'\in\cA'}I(a')=\sum_{p/2<q\le p}1\gg p\log^{-1}p.
$$
Therefore, there exists $a'_0\in\cA'$ such that
\begin{equation}
\label{eqn:I(a)} I(a'_0)\gg p|\cA'|^{-1}\log^{-1}p.
\end{equation}
Next, denote $\cA=\cA'\setminus\{a'_0\}$. Since $|\cA'|\ge 2,$ then
\begin{equation}
\label{eqn:AA} 0.5|\cA'|\le |\cA|\le |\cA'|.
\end{equation}
Now define the set
$$
\cB=\{\tau(q^2)\pmod p: \quad p/2<q\le p, \,\, \tau(q)\equiv
a'_0\pmod p\}.
$$
The elements of the set $\cB$ are of the form
$$
\tau(q^2)= \tau^2(q)-q^{11}\equiv {a'_0}^2-q^{11}\pmod p
$$
where $p/2<q\le p$ and $\tau(q)\equiv a'_0\pmod p.$ Thus, according
to~\eqref{eqn:I(a)} and~\eqref{eqn:AA}, $q$ runs through the set of
$$
\gg p|\cA'|^{-1}\log^{-1}p \gg p|\cA|^{-1}\log^{-1}p
$$
different residue classes modulo $p.$ Since $q^{11}\pmod p$ can take
any value at most 11 times, then
\begin{equation}
\label{eqn:|A||B|} |\cB|\gg p|\cA|^{-1}\log^{-1}p.
\end{equation}

Fix $\varepsilon$, $0<\varepsilon<0.1$. Let now $\cC$ be the set of
all different elements of the sequence of residues
$$
\tau(q)\pmod p, \quad \tau(q^2)\pmod p,
$$
where $q\le p^{0.5\varepsilon}.$ The above argument applied to the
sets
$$
\{\tau(q)\pmod p: \quad q\le p^{0.5\varepsilon}\}
$$
and
$$
\{\tau(q^2)\pmod p: \quad q\le p^{0.5\varepsilon}\}
$$
shows that $|\cC|\gg p^{\varepsilon/6}.$  Note that the sets $\cA,
\cB, \cC$ are formed with elements of types $\tau(n_1),
\tau(n_2),\tau(n_3)$ correspondingly in such a way, that the numbers
$n_1, n_2, n_3$ are pairwise coprime and $n_1\le p, n_2\le p^{2},
n_3\le p^{\varepsilon}.$

\bigskip

If $|\cA|<p^{0.1\varepsilon},$ then by~\eqref{eqn:|A||B|}, $|\cB|\gg
p^{1-\varepsilon/9}$ and thus $|\cB||\cC|\gg p^{1+0.01\varepsilon}.$
Therefore we can apply Lemma~\ref{lem:Glib} with $\cX=\cB$ and
$\cY=\cC$ and use the multiplicative property of $\tau(n).$ The
elements of the set
$$
\cX\cY=\{xy: \quad x\in\cX, y\in\cY\}
$$
in this case will be of the form $\tau(n)\pmod p$ with $n\le
p^{2+\varepsilon}.$

\bigskip

If $|\cA|>p^{2/3},$ then we first split $\cA$ into two disjoint
subsets $\cA_1$ and $\cA_2$ with $\ge p^{2/3}/3$ elements in each.
Then apply Lemma~\ref{lem:Glib} with $\cX=\cA_1$ and $\cY=\cA_2.$
The elements of the set $\cX\cY$ in this case take the form
$\tau(n)\pmod p$ with $n\le p^{2}.$

\bigskip

If $p^{0.1\varepsilon}\le \cA\le p^{2/3},$ then denote by $\cT$ the
one of the sets $\cA+\cC$ and $\cA\cC$ with the biggest cardinality.
Since $|\cC|\gg p^{\varepsilon/6},$ according to Bourgain's
estimate~\cite[Theorem 1.1]{Bourg}, there exists a positive constant
$\gamma=\gamma(\varepsilon)>0$ such that
$$
|\cT|\gg p^{\gamma}|\cA|.
$$
Then $|\cB||\cT|\gg p^{1+\gamma/2}.$ Hence, we can apply
Lemma~\ref{lem:Glib} with $\cX=\cB$ and $\cY=\cT.$  The elements of
the set $\cX\cY$ in this case will be either of the form
$\tau(n_1)+\tau(n_2)\pmod p$ with $n\le p^{3}$ or of the form
$\tau(n)\pmod p$ with $n\le p^{3+\varepsilon}.$

\bigskip

The result now follows.


\begin{thebibliography}{9999}

\bibitem{ArCh} G. I. Arkhipov and V. N. Chubarikov, {\it O chisle slagayemix v additivnoy probleme Vinogradova
i ee obobsheniyax}, IV Mezhdunarodnaya konferentsiya ``Sovremenniye
problemy teorii chisel i ee prilozheniya", Aktualniye problemy,
Chast I, 5--38 (2001).

\bibitem{Bourg} J. Bourgain, {\it More on the sum-product phenomenon in prime fields and its
applications,} Int. J. Number Theory, {\bf 1} (2005), 1--32.


\bibitem{BKT} J. Bourgain, N. Katz and T. Tao, {\it  A sum-product estimate in finite fileds and
their applications,} Geom. Func. Anal., {\bf 14} (2004), 27--57.

\bibitem{Del} P. Deligne, {\it La conjecture de Weil} I, (French) Inst. Hautes Études Sci. Publ. Math.
No. {\bf 43} (1974), 273--307.

\bibitem{Gl} A. A. Glibichuk, {\it Combinational properties of sets of residues modulo a prime
and the Erd\H{o}s-Graham problem,} Mat. Zametki, {\bf 79} (2006),
384--395; translation in: Math. Notes, {\bf 79} (2006), 356--365.

\bibitem{Hua} L. K. Hua, `The additive prime number theory',
(Russian), 1947.

\bibitem{Iw} H. Iwaniec, {`Topics in Clasical Automorphic Forms',}
AMS, Providence, Rhode Island, 1997.

\bibitem{Kar} A. A. Karatsuba, {\it Additive congruences,} (Russian), Izv. Ross.
Akad. Nauk Ser. Mat. {\bf 61} (1997), no. 2, 81--94; translation in
Izv. Math., {\bf 61} (1997), no. 2, 317--329.

\bibitem{Kobl} N. Koblitz, ``Introduction to elliptic curves and modular
forms", Springer-Verlag, New York, 1993.

\bibitem{KT} A. V. Kumchev and D. I. Tolev, {\it An invitation to additive prime number
theory}, Serdica Math. J., {\bf 31} (2005), 1--74.

\bibitem{Nie} D. Niebur, {\it A formula for Ramanujan's $\tau
$-function,} Illinois J. Math., {\bf 19} (1975), 448--449.


\bibitem{Ser} J. P. Serre, {\it Congruences et formes modulaires} [d'apr\`es H. P. F.
Swinnerton-Dyer] (French), Lecture Notes in Math, {\bf 317},
319--338, Springer, Berlin, 1973.

\bibitem{Sh} I. E. Shparlinski, {\it On the value set of the Ramanujan
function,} Arch. Math., {\bf 85} (2005), 508--513.


\vspace{1cm}

Adresses of the authors:
\bigskip

M.~Z.~Garaev and V. C. Garcia:\\
Instituto de Matem{\'a}ticas UNAM, C.P. 58089, Morelia,
Michoac{\'a}n, M{\'e}xico.

\bigskip

S.~V.~Konyagin:\\
Dept. of Mechanics and Mathematics,  Moscow State University,
Moscow, 119992, Russia.

\vspace{1cm}

emails:

garaev@matmor.unam.mx

garci@matmor.unam.mx

konyagin@ok.ru


\end{thebibliography}
\end{document}